\documentclass{amsart}
\usepackage{amssymb}
\usepackage{geometry}

\theoremstyle{definition}

\theoremstyle{remark}

\numberwithin{equation}{section}

\input{tcilatex}
\begin{document}
\title{On the K-Ring of the Classifying Space of the Generalized Quaternion
Group}
\author{MEHMET\ KIRDAR\ AND\ SEV\.{I}LAY\ \"{O}ZDEM\.{I}R}
\address{Department of Mathematics, Faculty of Arts and Science, Nam\i k
Kemal University, Tekirda\u{g}, Turkey}
\email{mkirdar@nku.edu.tr, gilthoniel1986@gmail.com}
\date{13 December, 2012}

\begin{abstract}
We describe the K-ring of the classifying space of the generalized
quaternion group in terms of generators and the minimal set of relations. We
also compute the order of the main generator in the truncated rings.
\end{abstract}

\maketitle







\section{Introduction}

The K-ring of the classifying space $BQ_{2^{n}}$ of the generalized
quaternion group $Q_{2^{n}}$, $n\geq 3,$ is described in [2] and [4]. In
this note, we will make the description of these rings in a simpler way, by
a minimal set of relations on a minimal set of generators. We will also make
nice connections of these computations with the computations done for the
lens spaces.

In particular, we will compute the order of the main generator in truncated
rings, in a much shorter way than that is done in [4].

As usual, the description of these rings will be done from the
representation rings of the groups via the Atiyah-Segal Completion Theorem
(ASCT) which roughly says that the K-ring of the classifying space of a
group is the completion of the representation ring of this group at the
augmentation ideal of the representation ring.

We will also check the minimality of the relations we found through the
Atiyah-Hirzebruch Spectral Sequence (AHSS), as usual.

In connection with and parallel to this problem, the reader is advised to
look at the computations done for cyclic and dihedral groups. A quick survey
for the $K$-Rings of the classifying spaces of cyclic groups, i.e. lens
spaces, and of the classifying spaces of dihedral groups can be found in [3].

\section{Representations}

The quaternion group $Q_{2^{n}}$, where $n\geq 3,$ is generated by two
elements $x$ and $y$ with the relations $x^{2^{n-1}}=1,$ $x^{2^{n-2}}=y^{2}$
and $xyx=y.$ Note that $x$ generates a cyclic group of order $2^{n-1}$ and $%
y $ generates a cyclic group of order $4.$ Note also that $x^{2^{n-2}},$
which is equal to $y^{2},$ generate a cyclic group of order $2.$ We keep in
mind these natural group inclusions of the cyclic groups in $Q_{2^{n}}.$

There are four one dimensional irreducible complex representations of $%
Q_{2^{n}}.$ We will denote them by $1,\eta _{1},\eta _{2}$ and $\eta _{3}.$
They are explained by $1\times 1$ matrices in [4].

There are $2^{n-2}-1$ two dimensional irreducible complex representations of 
$Q_{2^{n}}.$ We will denote them by $d_{i}$ where $1\leq i\leq 2^{n-2}-1.$
Actually, $d_{i}$ make sense for any integer $i$ and this will be clarified
below.

Since all we need will be the relations that they can generate, we will not
describe these representations by matrices here. The descriptions of these
representations by $2\times 2$ matrices are given in [4].

Before presenting the relations let us set: $2^{n}=2m=4k.$ And we have this
convention from now on. Note that $k\geq 2$ and it is a power of 2 too.

Now, we will list all possible relations in the representation ring $%
R(Q_{4k})$. First of all $\eta _{3}=\eta _{1}\eta _{2}.$ And since $\eta
_{1}^{2}=\eta _{2}^{2}=1$, we also have the relation $\eta _{3}^{2}=1$. For $%
d_{i}$'s, we have the start $d_{0}=1+\eta _{1}$ and we have the end $%
d_{k}=\eta _{2}+\eta _{3}.$

The main relation, which is the most important of all, is%
\begin{equation*}
d_{i}d_{j}=d_{i+j}+d_{i-j}.
\end{equation*}

This relation make sense for any integer couple $i,j$ because of the
following fact: $d_{i}=d_{m-i}$ for all integers $i.$

And another set of relations are for the products of the one dimensional and
two dimensional representations, and they are: $\eta _{1}d_{i}=d_{i}$ and $%
\eta _{2}d_{i}=d_{k-i}$ for all $i$. Since, $\eta _{3}=\eta _{1}\eta _{2},$
it follows that $\eta _{3}d_{i}=d_{k-i}$ for all $i$, the same as $\eta
_{2}. $

We deduce from the relations above that, the representation ring of $%
Q_{2^{n}}$ is just generated by $\eta _{1},\eta _{2}$ and $d_{1},$ by means
of tensor products and direct sums. And the minimal polynomials on $\eta
_{1},\eta _{2}$ and $d_{1}$ that define the ring can be found from these
relations. But, we will do that in our new variables.

\section{Cohomology}

Integral cohomology of $Q_{4k},$ $k\geq 2,$ is the following and can be
found in [1]: 
\begin{equation*}
H^{p}(BQ_{4k};Z)=\left\{ 
\begin{tabular}{l}
$Z$ \\ 
$Z_{2}\oplus Z_{2}$ \\ 
$Z_{4k}$ \\ 
$0$%
\end{tabular}%
\right. 
\begin{tabular}{l}
$p=0$ \\ 
if $p=4s+2$ \\ 
if $p=4s,$ $s\geq 1$ \\ 
if $p$ is odd.%
\end{tabular}%
\end{equation*}

Note that the odd dimensional cohomology vanishes. Because of that, the AHSS
which converges to $K(BQ_{4k}),$ collapses on page two so that the $K$-ring
is completely determined by the integral cohomology and vise versa. Here, we
also notice that the cohomology is periodic.

The relations of the cohomology ring can also be found in [1]. Note that
these relations are quite different than the relations in the $K$-ring. We
don't try to find connections between these relations. We will just compare
the orders of the elements of $K(BQ_{4k})$ in the filtrations of the
spectral sequence with the sizes of the cohomology groups to prove that the
relations are minimal.

\section{K-Rings}

Corresponding to the representations $\eta _{1},\eta _{2}$ and $d_{1},$
there are induced vector bundles over the classifying space $BQ_{4k}$ and we
denote them by the same letters. We will set the reduced vector bundles as $%
v_{1}=\eta _{1}-1,$ $v_{2}=\eta _{2}-1$ and, most importantly, the main
element $\phi =d_{1}-2.$ \ Due to the ASCT, the elements $v_{1},v_{2}$ and $%
\phi $ generate $K(BQ_{4k})$. All we need is to find the minimal relations
on these generators so that the ring is well described.

First off all, since $\eta _{1}^{2}=\eta _{2}^{2}=1,$ we have the following
relations 
\begin{equation*}
v_{1}^{2}=-2v\ \text{and }v_{2}^{2}=-2v_{2}\text{ (\textbf{Relations 1 \& 2}%
).}
\end{equation*}

We note that the above relations are the standard relations for real line
bundles over the classifying spaces. These small and simple relations on $%
v_{i}$'s explain the cohomology groups $H^{4s+2}(BQ_{4k};Z)=Z_{2}\oplus
Z_{2} $ which are the $(4s+2)$-th filtrations $E_{\infty }^{4s+2,-4s-2}$ on
the main diagonal of the AHSS. The first $Z_{2}$ in the direct sum is
generated by $v_{1}^{s}$ and the second $Z_{2}$ is generated by $v_{2}^{s}$
where $s\geq 1.$

Next we should explain $Z_{4k}$'s which occur in the cohomology ring, in
other words, we should explain the filtrations $E_{\infty }^{4s,-4s},$ $%
s\geq 1,$ of the AHSS. This will not be easy.

Recall the natural inclusion of $Z_{2k}$ in $Q_{4k}$ defined by the element $%
x\in Q_{4k}.$ This inclusion rises a natural ring homomorphism $%
K(BQ_{4k})\rightarrow K(BZ_{2k}).$

Under this homomorphism, the image of the virtual bundle $d_{i}-2$ in $%
K(BZ_{2k})$ is $\eta ^{i}+\eta ^{-i}-2$ where $\eta $ is the one and only
generator of $K(BZ_{2k})$ and the one and only relation it satisfies is $%
\eta ^{2k}=1.$ We set $w=\eta +\eta ^{-1}-2$ in $K(BZ_{2k})$. The element $w$
generates a subring of $K(BZ_{2k})$ which is isomorphic to the subring of $%
KO(BZ_{2k}),$ solely generated by $w.$ Actually it is almost isomorphic,
except a $Z_{2}$ direct summand which is generated by the tautological one
dimensional reduced real bundle, traditionally denoted by $\lambda .$ Note
that $KO(BZ_{2k})$ is the real topological K-theory of the space $BZ_{2k}.$

Hence, under the natural homomorphism mentioned above, the image of $d_{i}-2$
in $K(BZ_{2k})$ is $\psi ^{i}(w)$ where $\psi ^{i}$ is the Adams operation
of degree $i.$ Lets recall the effect of (the real) Adams operation of
degree $i,$ on the main generator $w$ of $KO(BZ_{2k}):$%
\begin{equation*}
\psi ^{i}(w)=\dsum\limits_{j=1}^{i}\frac{\binom{i}{j}\binom{i+j-1}{j}}{%
\binom{2j-1}{j}}w^{j}
\end{equation*}

We name the above polynomial "quadratic binomial of degree $i$" because of
its connection to the real part of a root of unity and because its
coefficients are polynomials of $i^{2}$. In particular, under the above ring
homomorphism, $\phi $ maps exactly on $w.$

On the other hand, in the ring $K(BQ_{4k}),$ we have $d_{k+1}-d_{k-1}=0.$
Practically, we observe that this gives a polynomial in $\phi $ of degree $%
k+1$, no one dimensional bundles involved. This is true for any $d_{i}$
where $i$ is odd. They can be written as a polynomial of $\phi $ and only $%
\phi .$

In the ring $KO(BZ_{2k}),$ we know that the main relation is $\psi
^{k+1}(w)-\psi ^{k-1}(w)=0,$ [3]. And, we deduce that $\psi ^{k+1}(\phi
)-\psi ^{k-1}(\phi )=0$ in the ring $K(BQ_{4k})$ too. We also conclude that 
\begin{equation*}
d_{i}-2=\psi ^{i}(\phi )\text{ when }i\text{ is odd.}
\end{equation*}
We can also prove this very important fact from the relations of the
representation ring without referring to lens spaces. But, this would take
longer. Lens spaces make this tricky.

The polynomials $g_{2k}(\phi )=\psi ^{k+1}(\phi )-$ $\psi ^{k-1}(\phi )$ are
given by the following series:%
\begin{equation*}
g_{2k}(\phi )=4k\phi +\sum_{j=2}^{k}\frac{2k^{2}+j-1}{(j-1)(2j-1)}\binom{%
k+j-2}{2j-3}\phi ^{j}+\phi ^{k+1}
\end{equation*}

Hence, the following relation is satisfied in $K(BQ_{4k})$: 
\begin{equation*}
g_{2k}(\phi )=0\text{ (\textbf{Relation 3}).}
\end{equation*}%
From this relation, we deduce that $\phi $ satisfies a relation in the form 
\begin{equation*}
4k\phi =f(\phi )\phi ^{2}
\end{equation*}%
where $f(\phi )\in K(BQ_{4k})$ is a virtual bundle generated by $\phi .$
This explains the fact that the $4$-th filtration $E_{\infty }^{4,-4}$ on
the last page of the main diagonal of AHSS is generated by maybe $\phi $ (or
maybe $\phi -v_{1}-v_{2}$ etc.) and is isomorphic to $%
H^{4}(BQ_{4k};Z)=Z_{4k}.$ We don't want to speculate much about the spectral
sequence.

And by multiplying this relation by powers of $\phi ,$ all $4s$-th
filtrations on the main diagonal of AHSS, i.e. all groups $%
H^{4s}(BQ_{4k};Z)=Z_{4k}$ in the cohomology, are similarly explained.

But, we are still not done! It turns out that the relation 3 is not minimal.
We will prove it when we talk about the minimal relation for the cross
product $v_{1}v_{2}.$ We also didn't still explain the products $v_{1}\phi $
and $v_{1}\phi $. Without these relations, the ring can not be completely
described, although more or less the filtrations of the diagonal of the AHSS
are explained.

Lets first find the minimal relations for the products $v_{i}\phi $ where $%
i=1$ or $2,$ and explain why they are not needed to occupy any place on the
AHSS.

From the relations $\eta _{1}d_{1}=d_{1},$ it immediately follows that%
\begin{equation*}
v_{1}\phi =-2v_{1}\text{ (\textbf{Relation 4}).}
\end{equation*}%
This takes care of the product $v_{1}\phi .$ Next we will take care of the
product $v_{2}\phi .$ From the relation, $\eta _{2}d_{1}=d_{k-1,}$ since $%
d_{k-1}=\psi ^{k-1}(\phi )+2,$ we obtain%
\begin{equation*}
v_{2}\phi =\psi ^{k-1}(\phi )-\phi -2v_{2}\text{ (\textbf{Relation 5}).}
\end{equation*}

So $v_{i}\phi $ where $i=1$ or $2$ are dependent variables, and we don't
have to search a place on the AHSS for them.

Finally lets explain what remained, in other words, lets find the minimal
relation for the cross product $v_{1}v_{2}.$ It turns out that the main
relation is not the Relation 3, but that one. In fact, we will throw away
our favorite relation, the Relation 3, from the minimal set of relations.

We will separate the cases $n=3$ and $n\geq 4$, since it turns out that $%
Q_{8}$ is a little different than the bigger generalized quaternion groups.

For $k=2,$ from the relation $d_{1}^{2}=d_{2}+d_{0}=1+\eta _{1}+\eta
_{2}+\eta _{3},$ we have%
\begin{equation*}
v_{1}v_{2}=4\phi +\phi ^{2}-2v_{1}-2v_{2}\text{ (\textbf{Relation 6, for }}%
n=3\text{).}
\end{equation*}

This is the main relation for $K(BQ_{8}).$ If we multiply this equation by $%
\phi +2$, an amazing thing happens and we find the Relation 3. In other
words, we can kick off the Relation 3 from the minimal list of relations
that describes the ring.

For $n\geq 4,$ starting from the relation $d_{1}^{2}=d_{2}+d_{0},$ by using
the relation $d_{i}^{2}=d_{2i}+d_{0}$ repeatedly, one inside the other, we
can obtain the relation $d_{k}-d_{0}=\psi ^{k}(\phi ).$ And along the way,
amusingly, we obtain the polynomials $\psi ^{i}(\phi )$ where $i$ is a power
of 2, in terms of $\psi ^{\frac{i}{2}}(\phi ).$ And thus, we have%
\begin{equation*}
v_{1}v_{2}=\psi ^{k}(\phi )-2v_{2}\text{ (\textbf{Relation 6, for }}n\geq 4%
\text{).}
\end{equation*}

We will throw away the Relation 3 for $n\geq 4$ too. We multiply the
Relation 6 above by $\phi +2$ and then we use the Relations 4 and 5 properly
in the equation we obtained, and again amazingly find the Relation 3. This
also removes doubts from the obscure explanations we did above when we
derived this relation trickily from the lens spaces.

We sum up everything in:

\textbf{Theorem 1: }$K(BQ_{2^{n}})$ is generated by $v_{1},v_{2}$ and $\phi $
with the minimal set of relations (1),(2),(4),(5) and (6) above.

\section{Orders}

In [4], Proposition 5.1, the order of the element $\phi $ in the truncated
ring $R(Q_{4k})\diagup \phi ^{2}R(Q_{4k}),$ after a lot of work by huge
matrices, is found as $4k.$ He used these orders, to answer some geometric
problems, namely problems about immersion of spaces like \textit{\ }$%
S^{4N+3}\diagup Q_{4k}$ in real Euclidean spaces. But, from the relation $%
4k\phi =f(\phi )\phi ^{2}$ we found above, this is evident.

Similarly, we can find the order of $\phi $ in the truncated ring $%
R(Q_{4k})\diagup \phi ^{N+1}R(Q_{4k})$ by a careful counting. In the
Relation $3,$ we observe that the coefficient of $\phi ^{2}$ is $\frac{%
k(2k^{2}+1)}{3}$ whose primary 2 factor is $k=2^{n-2}.$ Since the
coefficient of $\phi $ is $4k=2^{n}$, the jump between them is $\frac{2^{n}}{%
2^{n-2}}=4.$ So the total count up to $\phi ^{N+1}$ must be $4k.4^{N}.$
Hence, we have  

\textbf{Corollary 2:  }The order of $\phi $ in \textbf{\ }$K(S^{4N+3}\diagup
Q_{4k})$ is $2^{n+2N}.$

\end{document}